# On the distribution of consecutive composite odd numbers and twin primes

**October 31, 2018**

**WOLF Marc**, https://orcid.org/0000-0002-6518-9882

**WOLF François**, https://orcid.org/0000-0002-3330-6087

**VILLEMIN François-Xavier**

Independent researchers;

marc.wolf3@wanadoo.fr, francois.wolf@dbmail.com, fxvillemin@orange.fr

**Abstract**:

We study pairs of consecutive odd numbers through a straightforward indexing. We focus in particular on twin primes and their distribution. With a counting argument, we calculate the limit of an alternating sum that is equal to 1 which means there are few twin primes. Finally, we show how to find the possible congruences for these prime numbers.

On the distribution of consecutive composite odd numbers and twin primes

# Contents

# On the distribution of consecutive composite odd numbers and twin primes

## Introduction

With the exception of the pair $(3,5)$, pairs of twin primes are of the form $(6m - 1, 6m + 1)$. This condition is obviously not sufficient to characterize them, and we propose to study here the distribution of pairs of composite numbers among them.

## I- Preliminaries

Let us recall here some of the notations from [1] and [2]:

1. $I$ is the set of odd integers greater than 1, i.e.:

$$I = \{N_k = 2k + 3; k \in \mathbb{N}\}$$

with $k$ the <u>index</u> of an odd number $N_k$. We note that $k$ is a multiple of 3 if and only if $N_k$ is;

2. $P$ is the set of odd prime numbers, primes will also be enumerated in ascending order $\wp_0 = 3, \wp_1, \ldots, \wp_n \ldots$, with $q_0 = 0, q_1, \ldots, q_n \ldots$ their respective indices – we also note $Q$ the set of indices of prime numbers. Apart from $\wp_0$, none of the $\wp_i$ is a multiple of 3, which means that for $i > 0$, $q_i$ is not a multiple of 3;

3. $C$ is the set of composite odd integers greater than 1, i.e.:

$$C = I \backslash P = \{N_k \in I \mid \exists (a,b) \in I, N_k = ab\}$$

4. The function $f: k \in \mathbb{N} \mapsto N_k \in I$ is bijective. The inverse function is $f^{-1}: N_k \in I \mapsto k = \frac{N_k - 3}{2}$. The inverse image of $C$ is called $W$:

$$W = f^{-1}(C) = \{k \in \mathbb{N} \mid N_k \in C\}$$

5. The function

$$k: (n, j) \in \mathbb{N}^* \times \mathbb{N} \mapsto k_j(n) = (2j + 3)n + j$$

is a surjection on $W$. In other words, $W$ is the (non-disjoint) union of the sets $W_j = \{k_j(n)\}_{n \in \mathbb{N}^*}$.

6. The remarkable indices are the indices of the form $k_j(j + 1) = 2j^2 + 6j + 3$, they correspond to the indices of odd squares. These indices are never of the form $3n + 1$. The counting interval $I_D(j)$ is defined as $I_D(j) = [\![0, k_{j+1}(j + 2)]\!]$. The cardinal of a set $E$ being noted $|E|$, we have $|I_D(j)| = k_{j+1}(j + 2) + 1 = 2j^2 + 10j + 12$.

7. We define $A(j)$ as the set of indices of $I_D(j)$ non-multiples of 3: $A(j) = \{j \in I_D(j) \mid j \bmod 3 \neq 0\}$. Likewise, we introduce $B(j) = A(j) \cap W$.

8. With $x \in \mathbb{R}$, we note $\lfloor x \rfloor$ its integer part, and $\lfloor x \rfloor_+$ its positive integer part. The number of integers in $\mathbb{N}$ not greater than $x$ is $\lfloor x + 1 \rfloor_+$.

To avoid confusion, we reserve uppercase Latin letters for sets of odd numbers or indices of these numbers. Sets of pairs of numbers will be referred to with uppercase Greek letters.

## II- Pair of indices of the form $(3n + 1, 3n + 2)$

### II.1 The set Γ and its remarkable subsets Π, Ψ and Ξ

We define $\Gamma(j)$ as the set of pairs of consecutive indices of $I_D(j)$ non-multiples of 3:

$$\Gamma(j) = \{(k_1, k_2) \in A(j)^2 \mid k_2 - k_1 = 1\}.$$

**Property 2.1**: $|A(j)| = 2 * |\Gamma(j)|$.

*Proof*: Since the remarkable index that is the largest element of $A(j)$ is never of the form $3n + 1$, we deduce that $3n + 1 \in A(j)$ if and only if $3n + 2 \in A(j)$. $A(j)$ has therefore as many numbers congruent to 1 as to 2 modulo 3. We note $A_1(j)$ (respectively $A_2(j)$) the numbers that are congruent to 1 (respectively 2) modulo 3.

It is thus clear that $|A(j)| = 2|A_1(j)| = 2|A_2(j)|$. Moreover, any element $\gamma \in \Gamma(j)$ is necessarily of the form $\gamma = (3n + 1, 3n + 2)$, so there is a bijection:

$$\begin{array}{ccc} A_1(j) & \to & \Gamma(j) \\ x & \mapsto & (x, x+1) \end{array}$$

which proves the property.

**Definition 2.1.2**: We define $\Pi(j)$, $\Psi(j)$ et $\Xi(j)$ as follows:

$$\Pi(j) = \Gamma(j) \cap (Q \times Q)$$
$$\Psi(j) = \Gamma(j) \cap (W \times W)$$
$$\Xi(j) = \Gamma(j) \setminus (\Pi(j) \cup \Psi(j))$$

These three sets clearly form a partition of $\Gamma(j)$. $\Pi(j)$ corresponds to pairs of twin primes less than $(2j + 5)^2$ with the exception of $(3,5)$. $\Psi(j)$ corresponds to the pairs of consecutive odd composite numbers, while $\Xi(j)$ corresponds to the mixed pairs.

## II.2 A counting result

We show here a property linking the cardinal of the different sets defined previously.

**Property 2.2**: We have the following equality:
$$|\Pi(j)| = |\Gamma(j)| - |B(j)| + |\Psi(j)|$$

*Proof*: The partition of $\Gamma(j)$ (See definition 2.1.2) allows us to write the following equality:
$$|\Gamma(j)| = |\Pi(j)| + |\Psi(j)| + |\Xi(j)|$$

Moreover, using the same arguments as in the proof of property 2.1, we have:
$$|B(j)| = 2 * |\Psi(j)| + |\Xi(j)|$$

We deduce that $|\Xi(j)| = |B(j)| - 2 * |\Psi(j)|$, and thus the desired result.

## III- Counting of $\Gamma(j)$, $B(j)$ and $\Psi(j)$

The previous property establishes a relation on the number of pairs of twin primes. We will now expand on the cardinals of the sets involved in this relation.

## III.1 Counting of $\Gamma(j)$ and $B(j)$

The property 3.1 of [2] gives us that for all $j \in \mathbb{N}$:

$$|A(j)| = \left\lfloor \frac{2}{3}(k_{j+1}(j+2) + 1) \right\rfloor = \frac{4}{3}j^2 + \frac{20}{3}j + O(1)$$

In particular:

$$|A(j)| \sim \frac{4}{3}j^2.$$

The following result is an immediate consequence.

On the distribution of consecutive composite odd numbers and twin primes

**Property 3.1**:
$$|\Gamma(j)| = \left\lfloor \frac{1}{3}(k_{j+1}(j+2)+1) \right\rfloor = \frac{2}{3}j^2 + \frac{10}{3}j + O(1)$$

In particular:
$$|\Gamma(j)| \sim \frac{2}{3}j^2.$$

Let us also remind the results of [2] on the cardinal of $B(j)$. If we denote by $\pi(x)$ the number of primes less than $x$, $\pi''(j)$ the number of indices of odd primes less than $(2j+5)$ i.e. $\pi''(j) = \pi((2j+5)) - 1$, and $\wp_K = \prod_{n \in K} \wp_n$, the underline{property 3.2.3} of [2] expresses $|B(j)|$ in two different ways:

$$|B(j)| = \sum_{\substack{K \subset [\![0, \pi''(j)]\!] \\ K \neq \emptyset, K \neq \{0\}}} (-1)^{|K|-1} \left\lfloor \frac{(2j+5)^2 - \wp_K}{2\wp_K} \right\rfloor$$

$$|B(j)| = -\sum_{k=1}^{k_{j+1}(j+2)} \mu(2k+3) \left\lfloor \frac{k_{j+1}(j+2) - k}{2k+3} \right\rfloor$$

The underline{property 3.2.4} of [2] also gives the asymptotic expansion of $|B(j)|$:
$$|B(j)| = |A(j)| - \pi((2j+5)^2) + 2 = \frac{4}{3}j^2 - \frac{2j^2}{\ln(j)} + o\left(\frac{j^2}{\ln(j)}\right)$$

## III.2 Counting of $\Psi(j)$

### III.2.1 The inclusion-exclusion principle
Let $(x, x+1) \in \Psi(j) = \Gamma(j) \cap W^2$. Thus $2x + 3 \leq (2j+5)^2$, and $2x + 3$ is a composite number, so it admits at least one odd prime factor between 5 and $2j + 5$. It's the same for his successor $2x + 5$.

We deduce the existence of two indices of prime numbers $q_s$ and $q_t$, with $1 \leq s, t \leq \pi''(j)$, such as $x \in W_{q_s}$ and $x + 1 \in W_{q_t}$.

We can notice that $q_s$ and $q_t$ are necessarily distinct: in fact, two consecutive odd numbers are necessarily coprime, and therefore $2x + 3$ and $2x + 5$ cannot admit a common prime divisor.

Furthermore, we may involve more prime numbers without changing the result. We deduce that more generally:

$$\forall N \geq \pi''(j) \quad \Psi(j) = \bigcup_{\substack{1 \leq s, t \leq N \\ s \neq t}} \Gamma(j) \cap (W_{q_s} \times W_{q_t}) \quad (1)$$

The inclusion-exclusion principle allows us to calculate the cardinal of this union:

$$\forall N \geq \pi''(j) \quad |\Psi(j)| = \sum_{\substack{K \subset [\![1,N]\!]^2 \\ K \neq \emptyset}} (-1)^{|K|-1} \left| \left( \bigcap_{(k1,k2) \in K} (W_{q_{k1}} \times W_{q_{k2}}) \right) \cap \Gamma(j) \right|. \quad (2)$$

On the distribution of consecutive composite odd numbers and twin primes

This result suggests to calculate, for all $K \subset [\![1, N]\!]^2$, the cardinal of $\bigcap_{(k_1,k_2)\in K}(W_{q_{k_1}} \times W_{q_{k_2}}) \cap \Gamma(j)$. So let us consider a subset K of $[\![1, N]\!]^2$, and analyze the properties of a pair $(x, x+1) \in \bigcap_{(k_1,k_2)\in K}(W_{q_{k_1}} \times W_{q_{k_2}}) \cap \Gamma(j)$.

Necessarily:
- $\forall (k_1, k_2) \in K \quad p_{k_1} | 2x + 3$ and $p_{k_2} | 2x + 5$
- $2x + 3 = 2 \mod 3$
- $2x + 3 = 1 \mod 2$

Let $K_1 = \{k_1 | \exists k_2 \ (k_1, k_2) \in K\}$ and also $K_2 = \{k_2 | \exists k_1 \ (k_1, k_2) \in K\}$. If $k \in K_1 \cap K_2$, then from the first point we deduce that $p_k$ divides both $2x + 3$ and $2x + 5$ that is not possible, therefore necessarily $K_1 \cap K_2 = \emptyset$.

Provided that this condition is verified, we can rewrite the equations above as a system of Diophantine equations on $y = 2x + 3$:
- $\forall k \in K_1 \quad y = 0 \mod p_k$
- $\forall k \in K_2 \quad y = -2 \mod p_k$
- $y = 2 \mod 3$
- $y = 1 \mod 2$

Insofar as the above congruences involve distinct two-by-two prime numbers, we can deduce (according to the Chinese theorem) the existence of a unique solution $y_0$ between 1 and $6 p_{K_1} p_{K_2} - 1$, such that every other solution in $\mathbb{N}$, and especially $y$, is of the form $y_0 + 6 m p_{K_1} p_{K_2}$, for $m \in \mathbb{N}$.

Conversely, it is easy to check that the conditions are sufficient provided that $2x + 3$ and $2x + 5$ are composite numbers, so that two cases are distinguished:

<u>Case 1</u>: If $y_0$ or $y_0 + 2$ is prime, then the first pair $(y_0, y_0 + 2)$ is excluded and necessarily:

$$\bigcap_{(k1,k2)\in K} (W_{q_{k_1}} \times W_{q_{k_2}}) \cap \Gamma(j)$$
$$= \left\{ \left( \frac{y_0 - 3}{2} + 3 m p_{K_1} p_{K_2}, \frac{y_0 - 1}{2} + 3 m p_{K_1} p_{K_2} \right), m \in \mathbb{N}^* \right\} \cap I_D(j).$$

In this case, the first coordinate of the smallest element is: $x_0 = \frac{y_0 - 3}{2} + 3 p_{K_1} p_{K_2}$.

<u>Case 2</u>: Otherwise, all the solutions of the Diophantine equation are suitable:

$$\bigcap_{(k1,k2)\in K} (W_{q_{k_1}} \times W_{q_{k_2}}) \cap \Gamma(j)$$
$$= \left\{ \left( \frac{y_0 - 3}{2} + 3 m p_{K_1} p_{K_2}, \frac{y_0 - 1}{2} + 3 m p_{K_1} p_{K_2} \right), m \in \mathbb{N} \right\} \cap I_D(j).$$

In this case, the first coordinate of the smallest element is: $x_0 = \frac{y_0 - 3}{2}$.

The primality of $y_0$ is possible only if $y_0 = p_{K_1}$ which also implies that $K_1$ must be reduced to one element. In the same way, the primality of $y_0 + 2$ implies $y_0 + 2 = p_{K_2}$ and is only possible when $K_2$ is a singleton. These are not sufficient conditions, but imply that case 1 occurs only for relatively low values of $y_0$, that is to say of the same order of magnitude as the prime numbers chosen for counting.

$y_0$ only depends on $K_1$ and $K_2$. Subsequently, when necessary, this dependency will be explicitly shown using the notation $y_0(K_1, K_2)$; and similarly for $x_0$ the first coordinate of the smallest element of $\bigcap_{(k_1,k_2)\in K}(W_{q_{k_1}} \times W_{q_{k_2}}) \cap \Gamma(j)$:

$$x_0(K_1, K_2) = \begin{cases} \dfrac{y_0(K_1, K_2) - 3}{2} + 3\wp_{K_1}\wp_{K_2} & \text{in case 1,} \\ \dfrac{y_0(K_1, K_2) - 3}{2} & \text{otherwise.} \end{cases}$$

**Property 3.2.1**: We have the equality $y_0(K_2, K_1) = 6\wp_{K_1}\wp_{K_2} - y_0(K_1, K_2) - 2$. Therefore, case 1 cannot occur for both $(K_1, K_2)$ and $(K_2, K_1)$.

*Proof*: $y_0(K_2, K_1)$ is the only solution of a system of Diophantine equations between 0 and $6\wp_{K_1}\wp_{K_2} - 1$, just check that $6\wp_{K_1}\wp_{K_2} - y_0(K_1, K_2) - 2$ is the solution of these same equations. We have:

$$\begin{cases} \forall k \in K_1 & y_0(K_1, K_2) = 0 \bmod \wp_k \\ \forall k \in K_2 & y_0(K_1, K_2) = -2 \bmod \wp_k \\ y_0(K_1, K_2) = 2 \bmod 3 \\ y_0(K_1, K_2) = 1 \bmod 2 \end{cases}$$

hence we deduce:

$$\begin{cases} \forall k \in K_1 & 6\wp_{K_1}\wp_{K_2} - y_0(K_1, K_2) - 2 = -2 \bmod \wp_k \\ \forall k \in K_2 & 6\wp_{K_1}\wp_{K_2} - y_0(K_1, K_2) - 2 = 0 \bmod \wp_k \\ 6\wp_{K_1}\wp_{K_2} - y_0(K_1, K_2) - 2 = 2 \bmod 3 \\ 6\wp_{K_1}\wp_{K_2} - y_0(K_1, K_2) - 2 = 1 \bmod 2. \end{cases}$$

Moreover, $y_0(K_1, K_2)$, as a multiple of $\wp_{K_1}$, cannot be equal to $6\wp_{K_1}\wp_{K_2} - 1$, that is to say that $6\wp_{K_1}\wp_{K_2} - y_0(K_1, K_2) - 2$ is also between 0 and $6\wp_{K_1}\wp_{K_2} - 1$, and therefore coincides with $y_0(K_2, K_1)$.

Finally, assuming that $y_0(K_1, K_2)$ is prime, this implies that it is equal to $\wp_{K_1}$, and that $\wp_{K_2}$ divides $\wp_{K_1} + 2$. Therefore, using that $\wp_{K_1}$ and $\wp_{K_2} \geq 5$:

$$\begin{aligned} y_0(K_2, K_1) &= (6\wp_{K_2} - 1)\wp_{K_1} - 2 \\ &\geq 29\wp_{K_1} - 2 \\ &\geq 28\wp_{K_1} + 3 \\ &> \wp_{K_1} + 2 \\ &\geq \wp_{K_2}. \end{aligned}$$

$y_0(K_2, K_1)$, multiple strict of $\wp_{K_2}$, cannot be prime, and similarly his successor $y_0(K_2, K_1) + 2$ cannot either. Similar inequalities can also be written in the case where $y_0(K_1, K_2) + 2$ is prime.

### III.2.2 Calculation of the cardinal of $\Psi(j)$

**Property 3.2.2**: The cardinal of $\bigcap_{(k1,k2)\in K}(W_{q_{k_1}} \times W_{q_{k_2}}) \cap \Gamma(j)$ is:

$$z(K_1, K_2) = \left\lfloor \frac{k_{j+1}(j+2) - x_0(K_1, K_2)}{3\wp_{K_1}\wp_{K_2}} + 1 \right\rfloor_+ . \quad (3)$$

*Proof*: We know that $(x, x+1) \in \bigcap_{(k_1,k_2)\in K} \left(W_{q_{k_1}} \times W_{q_{k_2}}\right) \cap \Gamma(j)$ if and only if there is $m \in \mathbb{N}$ such that:
$$x = x_0(K_1, K_2) + 3m\wp_{K_1}\wp_{K_2} \text{ and } x_0(K_1, K_2) + 3m\wp_{K_1}\wp_{K_2} \leq k_{j+1}(j+2).$$
This correspondence is clearly bijective and $x_0(K_1, K_2) + 3m\wp_{K_1}\wp_{K_2} \leq k_{j+1}(j+2)$ if and only if $m \leq \frac{k_{j+1}(j+2) - x_0(K_1,K_2)}{3\wp_{K_1}\wp_{K_2}}$. The result follows from point 8 of the preliminaries.

<u>Remark</u>: In the case 2, we know that $x_0(K_1, K_2) < 3\wp_{K_1}\wp_{K_2}$ and therefore:
$$\frac{k_{j+1}(j+2) - x_0(K_1, K_2)}{3\wp_{K_1}\wp_{K_2}} > -1.$$
In the case 1, we know that $3\wp_{K_1}\wp_{K_2} < x_0(K_1, K_2) \leq 3\wp_{K_1}\wp_{K_2} + \max(q_{K_1}, q_{K_2}) < 6\wp_{K_1}\wp_{K_2}$. Moreover, we know that one of values $(q_{K_1}, q_{K_2})$ is an index of a prime number $q$, and that the first coordinate of the pair of indices is increased by $q$. If we know that $q \leq k_{j+1}(j+2)$, then:
$$\frac{k_{j+1}(j+2) - x_0(K_1, K_2)}{3\wp_{K_1}\wp_{K_2}} > -1.$$
We deduce that in some cases we can use the integer part instead of the positive integer part:

**Corollary 3.2.2.1:** Let $N \geq \pi''(j)$. The cardinal $|\Psi(j)|$ is written:
$$|\Psi(j)| = \sum_{\substack{K \subset [\![1,N]\!]^2 \\ K \neq \emptyset}} (-1)^{|K|-1} \left\lfloor \frac{k_{j+1}(j+2) - x_0(K_1, K_2)}{3\wp_{K_1}\wp_{K_2}} + 1 \right\rfloor_+ \quad (4)$$

or:

$$|\Psi(j)| = \sum_{\substack{K_1 \subset [\![1,N]\!] \\ K_1 \neq \emptyset}} \sum_{\substack{K_2 \subset [\![1,N]\!] \\ K_2 \neq \emptyset, K_1 \cap K_2 = \emptyset}} (-1)^{|K_1|+|K_2|} \left\lfloor \frac{k_{j+1}(j+2) - x_0(K_1, K_2)}{3\wp_{K_1}\wp_{K_2}} + 1 \right\rfloor_+ \quad (5)$$

Moreover if $\wp_N \leq (2j+5)^2$, we can drop the positive parts in (4) and (5).

*Proof*: (4) stems directly from (2) and (3). To get (5), we group the terms corresponding to the same pair $(K_1, K_2)$. Then it is sufficient to show that:
$$\sum_{\substack{L \subset [\![1,N]\!]^2 \\ L_1 = K_1 \\ L_2 = K_2}} (-1)^{|L|-1} = (-1)^{|K_1|+|K_2|}.$$

Let us fix the set $K_2 = \{j_1 \ldots j_n\}$. Let us calculate first:
$$f(K_1) = \sum_{\substack{L \subset [\![1,N]\!]^2 \\ L_1 \subset K_1 \\ L_2 = K_2}} (-1)^{|L|-1}.$$

To build an index $L$ of the sum above uniquely, it suffices to choose for each element $j \in K_2$ a non-empty subset $L_1^j$ of $K_1$ so that $L_1^j = \{i \in K_1 | (i,j) \in L\}$. Thus:

$$f(K_1) = \sum_{\substack{L_1^{j_1} \subset K_1 \\ L_1^{j_1} \neq \emptyset}} \cdots \sum_{\substack{L_n^{j_n} \subset K_1 \\ L_1^{j_n} \neq \emptyset}} (-1)^{|L_1^{j_1}| + \cdots + |L_1^{j_n}| - 1} = -\left( \sum_{\substack{L_1 \subset K_1 \\ L_1 \neq \emptyset}} (-1)^{|L_1|} \right)^n.$$

Using the binomial expansion, we get:

$$\sum_{\substack{L_1 \subset K_1 \\ L_1 \neq \emptyset}} (-1)^{|L_1|} = \sum_{l=1}^{|K_1|} \binom{|K_1|}{l} (-1)^l = (1-1)^{|K_1|} - 1 = -1.$$

Thus $f(K_1) = (-1)^{|K_2|+1}$.

This yields the result for $|K_1| = 1$, as in this case the condition $L_1 \subset K_1$ is equivalent to $L_1 = K_1$. For any greater value of $|K_1|$, the principle of inclusion-exclusion yields:

$$\sum_{\substack{L \subset [\![1,N]\!]^2 \\ L_1 = K_1 \\ L_2 = K_2}} (-1)^{|L|-1} = \sum_{\substack{L_1 \subset |K_1| \\ L_1 \neq \emptyset}}^{|K_1|} (-1)^{|K_1|-|L_1|} f(K_1) = (-1)^{|K_1|+|K_2|-1} \sum_{l=1}^{|K_1|} \binom{|K_1|}{l} (-1)^l.$$

With one more binomial expansion, this leads to the expected result – the observation on positive parts is a simple consequence of the previous remark.

Remark: The sum in (5) contains much fewer terms than in (4): the former has the order of $2^{2n}$ against $2^{n^2}$ for the latter! Grouping terms $K_1$ and $K_2$ of equal size, it can be rewritten as an double alternating sum, while (4) becomes a (simple) alternating sum.

The corollary below focuses on the case where the condition $N \geq \pi''(j)$ is dropped:
**Corollary 3.2.2.2:** Let $N \in \mathbb{N}$. We have the following inequality:
$$e_N(j) \leq |\Psi(j)| \quad (6)$$

with:

$$e_N(j) = \sum_{\substack{K \subset [\![1,N]\!]^2 \\ K \neq \emptyset}} (-1)^{|K|-1} \left[ \frac{k_{j+1}(j+2) - x_0(K_1, K_2)}{3 \wp_{K_1} \wp_{K_2}} + 1 \right]_+$$

$$= \sum_{\substack{K_1 \subset [\![1,N]\!] \\ K_1 \neq \emptyset}} \sum_{\substack{K_2 \subset [\![1,N]\!] \\ K_2 \neq \emptyset, K_1 \cap K_2 = \emptyset}} (-1)^{|K_1|+|K_2|} \left[ \frac{k_{j+1}(j+2) - x_0(K_1, K_2)}{3 \wp_{K_1} \wp_{K_2}} + 1 \right]_+ \quad (7)$$

Moreover, $e_N(j)$ is increasing in $N$ and in $j$.
*Proof*: This is a generalization of the previous results, counting $\Psi_N(j) = \bigcup_{\substack{1 \leq s,t \leq N \\ s \neq t}} \Gamma(j) \cap (W_{q_s} \times W_{q_t})$ instead of $\Psi(j)$.

Clearly $\Psi_N(j) \subset \Psi_{N+1}(j)$ as more elements lay in the latter union and clearly as well $\Psi_N(j) \subset \Psi_N(j+1)$ as $\Gamma(j)$ is increasing, hence $e_N(j)$ is increasing in both $N$ and $j$.
Remark: Similarly, if we define the set $B_N(j)$ of multiples of at least one of $\wp_i$, $i = 1 \ldots N$, we get:

On the distribution of consecutive composite odd numbers and twin primes

$$|B_N(j)| = \sum_{\substack{K \subset [\![0,N]\!] \\ K \neq \emptyset, K \neq \{0\}}} (-1)^{|K|-1} \left\lfloor \frac{(2j+5)^2 - \wp_K}{2\wp_K} \right\rfloor$$

We will reuse sets $\Psi_N(j)$ and $B_N(j)$ in the last part of this article.

### III.2.3 Asymptotic expansion of $|\Psi(j)|$

As we have seen in [2], the prime number theorem implies that:

$$\pi'(x) = \pi((2j+5)^2) - 2 \sim \frac{2j^2}{\ln(j)}.$$

Moreover, as we have seen that two pairs of twin primes (greater than 3) cannot have a common value:

$$|\Pi(j)| \leq \frac{\pi'(x)}{2} \sim \frac{j^2}{\ln(j)}.$$

The twin prime numbers infinity remains a conjecture to this day, in particular we do not know an equivalent to $|\Pi(j)|$.

**Property 3.2.3**: We have the following asymptotic expansion:
$$|\Psi(j)| = \frac{2}{3}j^2 + O\left(\frac{j^2}{\ln(j)}\right)$$

*Proof*: According to property 2.2, we have $|\Psi(j)| = |B(j)| - |\Gamma(j)| + |\Pi(j)|$ and so $|\Psi(j)| = |B(j)| - |\Gamma(j)| + O\left(\frac{j^2}{\ln(j)}\right)$. The asymptotic expansions of $|B(j)|$ and $|\Gamma(j)|$ allows us to conclude.

### III.3 A special alternate sum in the equivalent of $|C_j|$

In this section we focus on another equivalent of $|\Psi(j)|$.
A naive manipulation of the formula:

$$|\Psi(j)| = \sum_{\substack{K_1 \subset [\![1,\pi''(j)]\!] \\ K_1 \neq \emptyset}} \sum_{\substack{K_2 \subset [\![1,\pi''(j)]\!] \\ K_2 \neq \emptyset, K_1 \cap K_2 = \emptyset}} (-1)^{|K_1|+|K_2|} \left\lfloor \frac{k_{j+1}(j+2) - x_0(K_1, K_2)}{3\wp_{K_1}\wp_{K_2}} + 1 \right\rfloor_+$$

consists in summing the equivalents of each term of the sum, which gives:

$$|\Psi(j)| \sim \frac{2j^2}{3} \sum_{\substack{K_1 \subset [\![1,\pi''(j)]\!] \\ K_1 \neq \emptyset}} \sum_{\substack{K_2 \subset [\![1,\pi''(j)]\!] \\ K_2 \neq \emptyset, K_1 \cap K_2 = \emptyset}} \frac{(-1)^{|K_1|+|K_2|}}{\wp_{K_1}\wp_{K_2}}$$

However, we must be careful that, as the sum has more and more terms, this approach is not mathematically valid. However, corollary 3.2.2.2 allows us to manipulate finite sums.

**Property 3.3.1**: Let $E_N = \sum_{\substack{K_1 \subset [\![1,N]\!] \\ K_1 \neq \emptyset}} \sum_{\substack{K_2 \subset [\![1,N]\!] \\ K_2 \neq \emptyset \\ K_1 \cap K_2 = \emptyset}} \frac{(-1)^{|K_1|+|K_2|}}{\wp_{K_1}\wp_{K_2}}$. The sum $(E_N)$ is increasing and bounded by 1.

*Proof*: (7) allows us to assert that for $j \to +\infty$, $e_N(j) \sim \frac{2j^2}{3} E_N$. But for all $N$, $e_N(j) \leq e_{N+1}(j)$ hence, dividing by $\frac{2j^2}{3}$ and taking the limit, $E_N \leq E_{N+1}$, that is, the sequence $(E_N)$ is increasing.

Moreover, $e_N(j) \leq |\Psi(j)|$ so by using the asymptotic expansion of $|\Psi(j)|$ we also deduce the inequality $E_N \leq 1$.

**Property 3.3.2**: The alternating sum $(E_N)$ converges to 1:

$$E_N = \sum_{\substack{K_1 \subset [\![1,N]\!] \\ K_1 \neq \emptyset}} \sum_{\substack{K_2 \subset [\![1,N]\!] \\ K_2 \neq \emptyset, K_1 \cap K_2 = \emptyset}} \frac{(-1)^{|K_1|+|K_2|}}{\mathcal{P}_{K_1} \mathcal{P}_{K_2}} \xrightarrow{N \to \infty} 1$$

*Proof*: The previous property already shows that $(E_N)$ converges to a limit not greater than 1. To obtain the desired result, we observe that the terms of the sum depend only on the disjoint union $K_1 \cup K_2 \coloneqq K$ and that terms can be grouped accordingly:

$$\sum_{\substack{K_1 \subset [\![1,N]\!] \\ K_1 \neq \emptyset}} \sum_{\substack{K_2 \subset [\![1,N]\!] \\ K_2 \neq \emptyset, K_1 \cap K_2 = \emptyset}} \frac{(-1)^{|K_1|+|K_2|}}{\mathcal{P}_{K_1} \mathcal{P}_{K_2}} = \sum_{\substack{K \subset [\![1,N]\!] \\ K \neq \emptyset}} (-1)^{|K|} \frac{1}{\mathcal{P}_K} \sum_{\substack{K_1 \subset K \\ K_1 \neq \emptyset, K_1 \neq K}} 1$$

$$= \sum_{\substack{K \subset [\![1,N]\!] \\ K \neq \emptyset}} (-1)^{|K|} \frac{2^{|K|} - 2}{\mathcal{P}_K}$$

$$= \sum_{\substack{K \subset [\![1,N]\!] \\ K \neq \emptyset}} \frac{(-2)^{|K|}}{\mathcal{P}_K} - 2 \sum_{\substack{K \subset [\![1,N]\!] \\ K \neq \emptyset}} \frac{(-1)^{|K|}}{\mathcal{P}_K}$$

$$= \left( \prod_{k=1}^{N} \left(1 - \frac{2}{\mathcal{P}_K}\right) - 1 \right) - 2 \left( \prod_{k=1}^{N} \left(1 - \frac{1}{\mathcal{P}_K}\right) - 1 \right)$$

$$\xrightarrow{N \to \infty} -1 + 2 = 1$$

In the equations above, we used the convergence of the infinite products $\prod_{k=1}^{N} \left(1 - \frac{2}{\mathcal{P}_K}\right)$ and $\prod_{k=1}^{N} \left(1 - \frac{1}{\mathcal{P}_K}\right)$ towards 0, which is a consequence of the divergence of the sum of the reciprocals of the primes (See [3]).

### III.4 The Möbius approximation

In the spirit of [2], we have obtained a Euler approximation of the proportion of pairs of composite odd numbers among the pairs of the form $(6m - 1, 6m + 1)$, namely $E_{\pi''(j)}$. It is also legitimate to study the Möbius approximation:

$$M_n = \sum_{k=1}^{n} \frac{\nu(2k+1)}{2k+1} \text{ with } \nu(n) = \begin{cases} (-1)^{|K|} [2^{|K|} - 2] & \text{if } n = \mathcal{P}_K \text{ with } 0 \notin K \\ 0 & \text{otherwise.} \end{cases}$$

An empirical study of this sum suggests its convergence, as shown in the graph below, made for $n$ ranging from 1 to 20 000.

On the distribution of consecutive composite odd numbers and twin primes

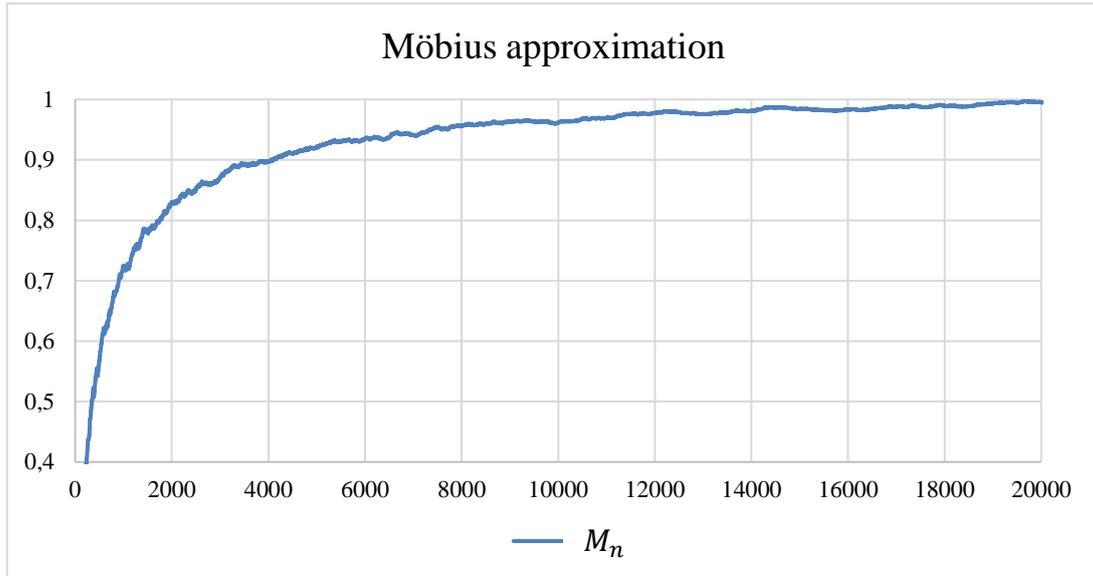

Möbius approximation

In addition to the Möbius function (see [4]), let us introduce the following three arithmetic functions, with $n = p_1^{\alpha_1} \ldots p_k^{\alpha_k}$:

$$\omega(n) = k; \quad \Omega(n) = \alpha_1 + \cdots + \alpha_k; \quad \mu_A(n) = A^{\omega(n)}\mu(n).$$

Thus, $\nu$ coincides with $\mu_2 - \mu$ on odd numbers non-multiple of 3.

The proof of [2] relating to Möbius approximation is difficult to adapt because it is based on a convergence result equivalent to the prime number theorem, and which requires analysis tools that are outside the scope of this article.

We can nevertheless show some interesting results on generalized the Möbius function $\mu_A$:

**Property 3.4.1**: The following identities are true for all $n \in \mathbb{N}^*$:

$$\sum_{d|n} \mu_A(d) = (1-A)^{\omega(n)}$$

$$\sum_{d|n} \mu_A(d).A^{\Omega(\frac{n}{d})} = \begin{cases} 1 \text{ if } n = 0, \\ 0 \text{ otherwise.} \end{cases}$$

*Proof*: Let's consider $n = p_1^{\alpha_1} \ldots p_{\omega(n)}^{\alpha_{\omega(n)}}$. Then:

$$\sum_{d|n} \mu_A(d) = \sum_{K \subset [\![1,\omega(n)]\!]} \mu_A(p_K) = \sum_{k=0}^{\omega(n)} \binom{n}{k}(-A)^k = (1-A)^{\omega(n)}.$$

The second equality is obvious for $n = 1$. Suppose first that n is the strict power of a prime number: $n = p^\alpha$. Then:

$$\sum_{d|n} \mu_A(d) A^{\Omega(\frac{n}{d})} = \sum_{\beta=0}^{\alpha} \mu_A(p^\beta) A^{\alpha-\beta} = A^\alpha - A.A^{\alpha-1} = 0.$$

Let's go back to the general case $n = p_1^{\alpha_1} \ldots p_{\omega(n)}^{\alpha_{\omega(n)}} = p_1^{\alpha_1} n'$. We observe that:

On the distribution of consecutive composite odd numbers and twin primes

$$\sum_{d|n} \mu_A(d) A^{\Omega(\frac{n}{d})} = \sum_{\substack{d_1|p_1^{\alpha_1} \\ d'|n'}} \mu_A(d_1 d') A^{\Omega\left(\frac{n}{d_1 d'}\right)}$$

$$= \sum_{\substack{d_1|p_1^{\alpha_1} \\ d'|n'}} \mu_A(d_1) \mu_A(d') A^{\Omega\left(\frac{p_1^{\alpha_1}}{d_1}\right) + \Omega\left(\frac{n'}{d'}\right)}$$

$$= \left(\sum_{d_1|p_1^{\alpha_1}} \mu_A(d) A^{\Omega(\frac{n}{d})}\right) \left(\sum_{d'|n'} \mu_A(d) A^{\Omega(\frac{n}{d})}\right)$$

$$= 0,$$

which allows us to conclude. We immediately deduce:

**Property 3.4.2**: For all $n \in \mathbb{N}^*$:

$$\sum_{1 \le k \le n} \mu_A(k) \left\lfloor \frac{n}{k} \right\rfloor = \sum_{1 \le j \le n} (1 - A)^{\omega(j)}.$$

*Proof*: It is a simple manipulation of sums:

$$\sum_{1 \le k \le n} \mu_A(k) \left\lfloor \frac{n}{k} \right\rfloor = \sum_{1 \le k \le n} \mu_A(k) \sum_{1 \le l \le \frac{n}{k}} 1$$

$$= \sum_{\substack{k,l \\ 1 \le kl \le n}} \mu_A(k)$$

$$= \sum_{1 \le j \le n} \sum_{k|j} \mu_A(k)$$

$$= \sum_{1 \le j \le n} (1 - A)^{\omega(j)}.$$

In the penultimate equality, we made a change of index $j = kl$ and in the last equality we used the property 3.4.1.

**Corollary 3.4.2**: $S_n := \sum_{k=1}^{n} \frac{\mu_2(k)}{k}$ is bounded by 2.

*Proof*: According to property 3.4.2 applied for $A = 2$, we have:

$$\left| \sum_{1 \le k \le n} \mu_2(k) \left\lfloor \frac{n}{k} \right\rfloor \right| = \left| \sum_{1 \le j \le n} (-1)^{\omega(j)} \right| \le n.$$

Moreover, for any $k$ between 1 and $n$:

$$\left\lfloor \frac{n}{k} \right\rfloor \le \frac{n}{k} < \left\lfloor \frac{n}{k} \right\rfloor + 1.$$

We then deduce:

$$\left| n \sum_{1 \leq k \leq n} \frac{\mu_2(k)}{k} \right| \leq \left| \sum_{1 \leq k \leq n} \mu_2(k) \left\lfloor \frac{n}{k} \right\rfloor \right| + n \leq 2n.$$

This allows us to conclude.

<u>Remark</u>: This corollary can be generalized for $A \in [0,2]$. For $A < 0$, $\mu_A$ is always positive, and it is easy to observe that $\sum_{1 \leq k \leq n} \frac{\mu_A(k)}{k}$ diverges knowing that the sum of reciprocals of prime numbers diverges.

We end with a last generalization of the result that shows that $M_n$ is bounded:

**Property 3.4.3**: For all $b \in \mathbb{N}^*$:

$$S_n^b := \sum_{\substack{1 \leq k \leq n \\ k \wedge b = 1}} \frac{\mu_2(k)}{k}$$

is bounded. Moreover, $M_n$ is bounded.

<u>Proof</u>: Let $\mu_{2,b}(n) = \mu_2(n)$ if $n \wedge b = 1$, 0 otherwise. We can generalize:

$$\sum_{d|n} \mu_{2,b}(d) = (-1)^{\omega_b(n)}$$

where $\omega_b(n)$ represents the number of prime divisors of $n$ not dividing $b$. We then deduce that $\sum_{1 \leq k \leq n} \mu_{2,b}(k) \left\lfloor \frac{n}{k} \right\rfloor = \sum_{1 \leq j \leq n} (-1)^{\omega_b(j)}$, and in particular $S_n^b$ is bounded. By fixing $b = 6$, it is easy to deduce that $M_n$ is bounded.

We thus leave the convergence of $M_n$ in the state of conjecture, and more generally that of the (bounded) series of the form, for $A \in [0,2]$:

$$\sum_{n \wedge b = 1} \frac{\mu_A(n)}{n}$$

## III.5 Fifteen possible congruences for twin primes

We have already seen that any prime number is written $6m - 1$ or $6m + 1$, with the exception of 2 and 3. The pairs of twin primes are therefore necessarily of the form $(6m - 1, 6m + 1)$, apart from $(3, 5)$.

Going further in congruences, we will generalize and demonstrate that there are 15 possible congruences modulo 105 for pairs of twin primes different from $(3, 5)$ and $(5, 7)$ pairs. None, of course, guarantees the primality of any of the two numbers.

Indeed, consider a pair of twin primes $(p, p + 2)$ such that $p$ is different from 3 and 5 (and therefore from 2 and 7 since $p + 2$ is prime). In particular, $p \geq 11$ and so $p$ and $p + 2$ are prime with $2, 3, 5, 7$ - that is, they are prime with $210 = 2 \times 3 \times 5 \times 7$.

A sieve method quickly eliminates all multiples of 2, 3, 5 and 7 between 0 and 209:

On the distribution of consecutive composite odd numbers and twin primes

| 0 | 1 | 2 | 3 | 4 | 5 | 6 | 7 | 8 | 9 | 10 | 11 | 12 | 13 | 14 |
|---|---|---|---|---|---|---|---|---|---|---|---|---|---|---|
| 15 | 16 | 17 | 18 | 19 | 20 | 21 | 22 | 23 | 24 | 25 | 26 | 27 | 28 | 29 |
| 30 | 31 | 32 | 33 | 34 | 35 | 36 | 37 | 38 | 39 | 40 | 41 | 42 | 43 | 44 |
| 45 | 46 | 47 | 48 | 49 | 50 | 51 | 52 | 53 | 54 | 55 | 56 | 57 | 58 | 59 |
| 60 | 61 | 62 | 63 | 64 | 65 | 66 | 67 | 68 | 69 | 70 | 71 | 72 | 73 | 74 |
| 75 | 76 | 77 | 78 | 79 | 80 | 81 | 82 | 83 | 84 | 85 | 86 | 87 | 88 | 89 |
| 90 | 91 | 92 | 93 | 94 | 95 | 96 | 97 | 98 | 99 | 100 | 101 | 102 | 103 | 104 |
| 105 | 106 | 107 | 108 | 109 | 110 | 111 | 112 | 113 | 114 | 115 | 116 | 117 | 118 | 119 |
| 120 | 121 | 122 | 123 | 124 | 125 | 126 | 127 | 128 | 129 | 130 | 131 | 132 | 133 | 134 |
| 135 | 136 | 137 | 138 | 139 | 140 | 141 | 142 | 143 | 144 | 145 | 146 | 147 | 148 | 149 |
| 150 | 151 | 152 | 153 | 154 | 155 | 156 | 157 | 158 | 159 | 160 | 161 | 162 | 163 | 164 |
| 165 | 166 | 167 | 168 | 169 | 170 | 171 | 172 | 173 | 174 | 175 | 176 | 177 | 178 | 179 |
| 180 | 181 | 182 | 183 | 184 | 185 | 186 | 187 | 188 | 189 | 190 | 191 | 192 | 193 | 194 |
| 195 | 196 | 197 | 198 | 199 | 200 | 201 | 202 | 203 | 204 | 205 | 206 | 207 | 208 | 209 |

This leaves 15 (highlighted) possibilities for congruence modulo 210 of $p$: 11, 17, 29, 41, 59, 71, 101, 107, 137, 149, 167, 179, 191, 197 and 209 (or -1). The central symmetry observed in the table above corresponds to the equality of proposition 3.2.1 which gives $y_0(K_2, K_1)$ as a function of $y_0(K_1, K_2)$.

The method is not specific to prime numbers chosen in the above: for any family of prime numbers $p_{k_1} \ldots p_{k_n}$, we can identify the congruences eligible for two twin prime numbers modulo $m \coloneqq 2 \times p_{k_1} \times \ldots \times p_{k_n}$. We observe that $-1$ and $1$ will always be prime with $m$, which ensures that there is always a solution. The symmetry of the table above is explained by the fact that, if $x$ and $x + 2$ are prime with $m$, the same is true for $m - x - 2$ and $m - x$. In particular, since $(-1,1)$ is the only fixed point of this symmetry, there is always an odd number of possible congruences.

We can simply prove that no congruence can guarantee the prime property of the two numbers: indeed, whatever $x$, $(m + 1)^2 x$ is congruent to $x$ modulo $m$, and is however clearly not prime !

The counting methods of the preceding sections also make it possible to predict the number of possible congruences. Indeed, consider $j$ such that $k_{j+1}(j + 2)$ is of the form $105n + 3$. A counting can be made using the notation of corollary 3.2.2.2 for $N = 2$, thus restricting to the multiplicity by 3, 5 and 7. The pairs obtained are therefore the "twin primes-candidates", ie the set $\Pi_2(j)$ of pairs whose two coordinates are prime with 2, 3, 5 and 7:

$$|\Pi_2(j)| = |\Gamma(j)| - |B_2(j)| + |\Psi_2(j)|.$$

According to property 3.1, we have:

$$|\Gamma(j)| = \left\lfloor \frac{1}{3}(k_{j+1}(j + 2) + 1) \right\rfloor = \left\lfloor \frac{105n + 4}{3} \right\rfloor = 35n + 1$$

$$|B_2(j)| = \sum_{\substack{K \subset [\![0,2]\!] \\ K \neq \emptyset, K \neq \{0\}}} (-1)^{|K|-1} \left\lfloor \frac{(2j+5)^2 - \mathcal{p}_K}{2\mathcal{p}_K} \right\rfloor$$

$$= \left\lfloor \frac{(2j+5)^2 - \mathcal{p}_1}{2\mathcal{p}_1} \right\rfloor + \left\lfloor \frac{(2j+5)^2 - \mathcal{p}_2}{2\mathcal{p}_2} \right\rfloor - \left\lfloor \frac{(2j+5)^2 - \mathcal{p}_0\mathcal{p}_1}{2\mathcal{p}_0\mathcal{p}_1} \right\rfloor$$

$$- \left\lfloor \frac{(2j+5)^2 - \mathcal{p}_0\mathcal{p}_2}{2\mathcal{p}_0\mathcal{p}_2} \right\rfloor - \left\lfloor \frac{(2j+5)^2 - \mathcal{p}_1\mathcal{p}_2}{2\mathcal{p}_1\mathcal{p}_2} \right\rfloor + \left\lfloor \frac{(2j+5)^2 - \mathcal{p}_0\mathcal{p}_1\mathcal{p}_2}{2\mathcal{p}_0\mathcal{p}_1\mathcal{p}_2} \right\rfloor$$

$$= 21n + 15n - 7n - 5n - 3n + n = 22n;$$

Finally:

$$|\Psi_2(j)| = \sum_{\substack{K_1 \subset [\![1,2]\!] \\ K_1 \neq \emptyset}} \sum_{\substack{K_2 \subset [\![1,2]\!] \\ K_2 \neq \emptyset, K_1 \cap K_2 = \emptyset}} (-1)^{|K_1|+|K_2|} \left\lfloor \frac{k_{j+1}(j+2) - x_0(K_1, K_2)}{3\mathcal{p}_{K_1}\mathcal{p}_{K_2}} + 1 \right\rfloor_+$$

$$= \left\lfloor \frac{105n + 3 - x_0(1,2)}{3\mathcal{p}_1\mathcal{p}_2} + 1 \right\rfloor + \left\lfloor \frac{105n + 3 - x_0(2,1)}{3\mathcal{p}_1\mathcal{p}_2} + 1 \right\rfloor = 2n.$$

For the last equation, the values of $x_0(1,2)$ and $x_0(2,1)$ calculated using the property 3.2.1 have been used.

Thus $|\Pi_2(j)| = 35n + 1 - 22n + 2n = 15n + 1$. This represents an isolated pair of twin primes, (5,7) of indices (1,2), and 15 congruences modulo 105 for the indices that correspond to the pairs previously found.

## Conclusion

Through simple counting of finite sets of indices of pairs of consecutive composite odd numbers, we have highlighted existence of an alternating sum convergent to 1, which reflects the fact that there are few twin primes. However, we leave open the question of the convergence of the Möbius approximation, as we could only prove its boundedness. Finally, we have developed a method which gives in particular 15 possibilities modulo 210 for the pairs of twin primes except $(3, 5)$ and $(5, 7)$.